\documentclass[12pt]{article}
\usepackage[cp866]{inputenc}
\usepackage{graphicx}
\usepackage{subfig}
\usepackage{extsizes}
\usepackage{epsfig}
\usepackage[numbers,sort&compress]{natbib}
\usepackage{citehack}
\usepackage{float}
\usepackage{amsmath}
\usepackage{amssymb}
\newcommand{\R}{\mathbb {R}}

\newtheorem{proposition}{Proposition}
\newtheorem{lemma}{Lemma}[section]

\makeatletter \@addtoreset{equation}{section} \makeatother
\date{\today}
\title{Estimates of some  integrals related to variations of smooth functions}
\author{Anatoly Neishtadt$^{1,2}$  \\
$^1$ Loughborough University, Loughborough, LE11 3TU, UK\\
 $^2$ Space Research Institute, Moscow, 113997, Russia}

\begin{document}
\maketitle

\begin{abstract}
Estimates of some integrals related to variations of smooth functions are presented.
\end{abstract}
\section{Statement of the problem and the result}

Let $Q$ be a closed cube in  $\R^n$, and $Q_1$ be a neighbourhood  of $Q$.   Consider a twice  differentiable function $\sigma\colon Q_1 \to \R\,$.  Suppose that  second  partial derivatives of $\sigma$ are  uniformly bounded, and   $|\sigma|<1$ in $ Q$.  Coordinates in  $\R^n$ are $x_1,\ldots ,x_n$.
Denote $d\mu=dx_1\ldots dx_n$ the volume element in  $\R^n$.

Let $f \colon (0,1]  \to \R\,$ be a strictly decreasing positive smooth function such that $f(\xi)\sqrt{\xi}$ is bounded and  $\int\limits_0^1\frac {f(\xi)}{\sqrt{\xi}}\,d\xi$ converges.

For $a>0$ denote
\begin{eqnarray}  
\label{sigma_a}
\widetilde{ {\partial \sigma} /{\partial x}}^a
=\begin{cases}
   {\partial \sigma}/{\partial x},\, \mbox{if}\ |\sigma|<a \,, \\
      0,\, \hskip 0.5cm \mbox{if}\ |\sigma | \ge a\,, \\
   \end{cases}   
 \end{eqnarray}
  \begin{eqnarray}  \label{gamma}
 \Gamma_1(a)=\int\limits_Q{\|\widetilde{ {\partial \sigma} /{\partial x}}^a\!\!\|}d\mu\, , \
 \Gamma_2(a)=\int\limits_Q\left[ \| {\partial \sigma}/{\partial x}\|
- \|  \widetilde{ {\partial \sigma}/{\partial x}}^a  \!\!\|\right]f(|\sigma|) d\mu\, . 
  \end{eqnarray}
   Integrals are understood in the sense of Riemann.  The integrand   in the last integral is defined by continuity at points  $x$ where $\sigma=0$: there it equals 0.        
\begin{proposition}           \label {L8.6}
For $0<a<1/2$
            $$
            \Gamma_1(a)<C_{1}\sqrt a, \quad  \Gamma_2(a)<C_{2}\,,
            $$
            where $C_{1}, C_{2}$ are positive constants, i.e. values that do not depend on $a$.
            \end{proposition} 
            This proposition with $f(\xi)=|\ln (\xi)|$ is contained in the thesis of the author \cite{k_diss}. It is presented here for the sake of reference. In  \cite{k_diss} it is used to obtain estimates of probability of capture into resonance in systems with degeneracies.

 \section{Proof of Proposition  \ref {L8.6} }
 In this proof $c_i$ are  positive constants, i.e. values that do not depend on $a$. Appearance such a constant in the proof means the assertion that  such a constant exists.
 
  \subsection{ Proof of existence of integrals}
Divide  $Q$ into cubes with a side  $l>0$. Denote    $Q_1$ the union  of cubes that contain points, where 
$\|\partial\sigma/\partial x\|\le l;  \ Q_2=Q\setminus Q_1$. Then  $\|\partial\sigma/\partial x\|<c_1l$ in  $Q_1$.  
  Because $\|\partial\sigma/\partial x\|>l$ in  $Q_2$,  the set 
 $\{x:\sigma(x)=a, \ x\in Q_2\}$ is a union of a finite number of  smooth hypersurfaces. Therefore the function  $\Phi=\|\widetilde{\partial\sigma/\partial x}^a\|$  is discontinuous in  $Q_2$ only on a 
 finite number of  smooth hypersurfaces. Thus this function  is integrable in $Q_2$. Therefore there exists a partition of $Q_2$ into sets such that  on each of these sets  the oscillation of   $\Phi$  does not exceed $l$. 
 Consider the union of this partition and $Q_1$. This is a partition of $Q$. On each set of this partition the oscillation of $\Phi$ does not exceed $l+c_1l$. Because  $l$  is arbitrarily small, $\Phi$ is integrable on  $Q$. Thus $\Gamma_1(a)$ exists. Similarly, one can prove that 
  $\Gamma_2(a)$ exists.
 
\subsection{ Reduction to the case of one variable} 
Let   $Q=\{(x_1,\ldots, x_{n})\ :\  |x_j-x_j^0|\le d, \ j=1,\ldots, n \}$. We have 
      \begin{eqnarray*}  \nonumber
 \Gamma_1(a)&=&\int\limits_{Q}{\|\widetilde{ \frac{\partial \sigma}{\partial x}}^a\!\!\|}d\mu=\int\limits_{Q}\sqrt{\sum\limits_{i=1}^{n}\left({|\widetilde{ \frac{\partial \sigma}{\partial x_i}}^a\!\!}|\right)^2}d\mu\le \\
 &\le&\int\limits_{Q}\left(\sum\limits_{i=1}^{n}{|\widetilde{ \frac{\partial \sigma}{\partial x_i}}^a\!\!}|\right)d\mu=\sum\limits_{i=1}^{n}\int\limits_{Q_i}(\prod\limits_{j=1,\, j\ne i}^{n}dx_j)\int\limits_{x_i^0-d}^{x_i^0+d}{|\widetilde{ \frac{\partial \sigma}{\partial x_i}}^a\!\!}|dx_i \,,
\end{eqnarray*}
where  $Q_i=\{(x_1,\ldots ,x_{i-1}, x_{i+1}, \ldots , x_{n})\ :\  |x_j-x_j^0|\le d, \ j\ne i\}$, and notation is similar to that in (\ref{sigma_a}).
Similarly, one can estimate  $\Gamma_2(a)$.
  \begin{lemma} \label {L8.6.1}
  Let a function  $\psi(z)$,  $0\le z\le d$, be twice differentiable and $ |\psi(z)|<1, \ |d^2\psi/dz^2|<c_1$. Denote for numbers $0\le\alpha<\beta\le1$ 
 $$
	^{\alpha}\widetilde{d\psi/dz}^{\beta}=
	\begin{cases} d\psi/dz,\ \mbox{\rm  if} \ |\psi(z)|\in[\alpha, \beta)\, , \\
	0, \ \hskip 0.5cm \mbox{\rm if} \ |\psi(z)|\notin[\alpha, \beta)
	\end{cases}
	$$
	(in particular, $ ^0\widetilde{d\psi/dz}^{\beta}=\widetilde{ {d\psi}/{dz}}^{\beta}, \ ^{\beta}\widetilde{ {d\psi}/{dz}}^1= {d\psi}/{dz}-\widetilde{ {d\psi}/{dz}}^{\beta}$)\,.
	Then
	  $$\int\limits_{0}^{d}{|\widetilde{ {d\psi}/{dz}}^a\!\!|}dz<c_2\sqrt a, \ \int\limits_{0}^{d}{|^a \widetilde{ {d\psi}/{dz}}^1\!\!|}f(|\psi |)dz<c_3\,,$$ where constants $c_2$ and   $c_3$ can be expressed via $c_1$ and $d$. 
	\end {lemma}
	This lemma and the transformation of integrals above imply Proposition \ref {L8.6}.

\subsection{ Proof of Lemma \ref{L8.6.1}} 
 $1^{\circ}$. Estimate  $\int\limits_{0}^{d}{|\widetilde{ {d\psi}/{dz}}^a\!\!|}dz$. To this end we will use in a modified form a construction from  \cite{anosov}.
   Let  $l>0$. Denote  $M$ the  set of points in $[0,d]$, such that  $\left |{d\psi}/{dz}\right|\le l$ in $M$. Denote  $\Lambda=[0,d]\setminus M$. 
   Introduce $h=0.5c_1^{-1}l$  and construct sets  $M_h$ and  $\Lambda_h$ as follows. Let 
  $
   0<h<...<qh<d\le(q+1)h\,.
   $
  Then $M_h$ is the union of those intervals  $[(r-1)h, rh], (r=1,...,q), [qh,d]$ that belong to  $M$. Denote
   $\Lambda_h=\overline{[0,d]\setminus M_h}$. Clearly,
   $$\int\limits_{M_h}{|\widetilde{ {d\psi}/{dz}}^a\!\!|}dz\le dl\,.$$
   \begin{lemma} \label {L8.6.1.1}  $\left |{d\psi}/{dz}\right|>0.5l^{-1}$  at points of $\Lambda_h$ .
\end {lemma}

\underline{Proof.} Let  $z'\in\Lambda_h$. Then there exists   $r$ such that  $z'\in[(r-1)h, rh]$ and not all points of this interval belong to $M$ (in the case  if  $r=q+1$, not all points of the interval  $[qh, d]$ belong to $M$). Then there exists  $z''\in\Lambda$ such that   $|z'-z''|<h$. Then
$$
|\left(\frac{d\psi}{dz}\right)_{z=z'}|>|\left(\frac{d\psi}{dz}\right)_{z=z''}|-c_1h>l-c_1 0.5c_1^{-1}l=0.5l
$$
as required.
\vskip 0.7cm

The set   $\Lambda_h$  can be represented as a union of intervals  $[\alpha_r, \beta_r]$,  $1\le r\le r_1\le q+1<dh^{-1}+1$. On  $[\alpha_r, \beta_r]$ the function  $\psi$ is monotonous  (because  $|{d\psi}/{dz}|>0.5l$). Therefore
  \begin{eqnarray*}  \nonumber
 \int\limits_{\alpha_r}^{\beta_r}{|\widetilde{ {d\psi}/{dz}}^a\!\!|}dz&=&| \int\limits_{\alpha_r}^{\beta_r}({\widetilde{ {d\psi}/{dz}}^a\!\!})dz| \le2a, \\
  \int\limits_{\Lambda_h}{|\widetilde{ {d\psi}/{dz}}^a\!\!|}dz&=&\sum\limits_{r=1}^{r_1}\int\limits_{\alpha_r}^{\beta_r}{|\widetilde{ {d\psi}/{dz}}^a\!\!|}dz \le2r_1a<2(dh^{-1}+1)a=4dc_1l^{-1}a+2a, \\
 \int\limits_{0}^{d}{|\widetilde{ {d\psi}/{dz}}^a\!\!|}dz&=& \int\limits_{M_h}{|\widetilde{ {d\psi}/{dz}}^a\!\!|}dz+ \int\limits_{\Lambda_h}{|\widetilde{ {d\psi}/{dz}}^a\!\!|}dz< d(l+4c_1l^{-1}a+2a)\,.   
\end{eqnarray*}
 Choose  $l=\sqrt a$. Then we get the required estimate
$$
 \int\limits_{0}^{d}{|\widetilde{ {d\psi}/{dz}}^a\!\!|}dz< d(l+4c_1)\sqrt a+2a<c_2\sqrt a\,.
 $$  
 \medskip
 
 $2^{\circ}$. Estimate  $\Gamma= \int\limits_{0}^{d}{|^{a}\widetilde{ {d\psi}/{dz}}^1\!\!|}f(|\psi|)dz$. Consider a partition   $\{\kappa_r\}$ of the interval  $[a, 1]$:
  $$a=\kappa_1<\kappa_2<\ldots<\kappa_{k+1}=1\,.
  $$
 Denote
 $$
 (\Delta\psi)_\nu= \int\limits_{0}^{d}{|^{\kappa_\nu}\widetilde{ {d\psi}/{dz}}^{\kappa_{\nu+1}}|}dz
 $$ 
  Then
   \begin{eqnarray*}  \nonumber   
 \Gamma=\sum\limits_{\nu=1}^{k} \int\limits_{0}^{d}{|^{\kappa_\nu}\widetilde{ {d\psi}/{dz}}^{\kappa_{\nu+1}}|}f(|\psi|)dz &\le& \sum\limits_{\nu=1}^{k}
  f({\kappa_\nu})(\Delta\psi)_\nu\,,\\
   \int\limits_{0}^{d}{|^{\kappa_1}\widetilde{ {d\psi}/{dz}}^{\kappa_{\nu+1}}|}dz&= &\sum\limits_{\nu=1}^r(\Delta\psi)_\nu\,.\\
\end{eqnarray*} 
  Result of   n. $1^{\circ}$ implies that
   $$  
   \int\limits_{0}^{d}{|^{\kappa_1}\widetilde{ {d\psi}/{dz}}^{\kappa_{r+1}}|}dz\le\int\limits_{0}^{d}{|^{0}\widetilde{ {d\psi}/{dz}}^{\kappa_{r+1}}|}dz <c_2\sqrt{\kappa_{r+1}}\,.  
 $$    
 Therefore $\sum\limits_{\nu=1}^r(\Delta\psi)_\nu<c_2\sqrt{\kappa_{r+1}}$\,.  

 Denote $y=(y_1,.., y_k), \ Y{(y)}=\sum\limits_{\nu=1}^k f(\kappa_\nu)y_\nu$,
 $$
  {\cal Y}=\{y:\sum\limits_{\nu=1}^ry_\nu\le c_2\sqrt{\kappa_{r+1}}, \ y_r\ge0 \, (r=1,\ldots, k)\}.
  $$
  \begin{lemma} \label {L8.6.1.2} 
  $\sup\limits_{y\in{\cal Y}}  Y{(y)}=   Y{(y^0)}$, where $y^0$  such that
  $\sum\limits_{\nu=1}^ry_\nu^0= c_2\sqrt{\kappa_{r+1}}, \ r=1,\ldots, k$ (or, which is the same,
  $
  y_1^0= c_2\sqrt{\kappa_2}, \  y_\nu^0= c_2(\sqrt{\kappa_{\nu+1}}-\sqrt{\kappa_\nu}), \  \nu=2,\ldots, k)\,.
  $
\end {lemma}
   \begin{lemma} \label {L8.6.1.3} Denote  $\theta$ the  diameter of  partition  $\{\kappa_r\}$.
  Then  $\lim\limits_{\theta\to0}Y{(y^0)}<c_3$.
\end {lemma}

Because  $\Gamma\le\sup\limits_{y\in{\cal Y}}Y{(y)}=Y{(y^0)}$ and  $\Gamma$ does not depend on a choice of partition  $\{\kappa_r\}$, in the limit as    ${\theta\to0}$ we get  $\Gamma<c_3$, as required.
\section{Proofs of Lemmas \ref{L8.6.1.2} and \ref{L8.6.1.3} }
\label{last}
 {\underline {Proof of Lemma  \ref{L8.6.1.2} .}} Denote $\chi_r=c_2\sqrt{\kappa_{r+1}}$. The considered supremum is attained at some point  $y^0\in{\cal Y}$. Suppose that  $y_1^0<\chi_1$. There exists an index  $\nu_1 \,\,(1<\nu_1\le k)$  such that  $y^0_{\nu_1}\ne 0$.  Let  $\nu_1=2$ for definiteness.  
 Denote  $\delta=\min\{\chi_1-y_1^0, y_2^0\}>0, \ \tilde y^0=(y_1^0+\delta,y_2^0-\delta, y_3^0,..., y_k^0)$ . Clearly  $\tilde y^0\in{\cal Y}$. Further,
 $$
Y(\tilde y^0)=\delta(f(\kappa_1)- f(\kappa_2))+Y(y^0)>Y(y^0),
$$
which contradicts to the definition of  $y^0$. Therefore
$$
y_1^0=\chi_1=c_2\sqrt{\kappa_2}\quad  \mbox{and}\quad  \sum\limits_{\nu=2}^ry^0_\nu\le\chi_r-\chi_1\quad (r=2,\ldots, k)
$$
Repeating this reasoning we find that 

$$
y^0_\nu=\chi_\nu-\chi_{\nu-1}=c_2(\sqrt{\kappa_{\nu+1}}-\sqrt{\kappa_\nu}), \ \nu=2,\ldots, k.
$$

 {\underline {Proof of Lemma   \ref{L8.6.1.3} .}}
  \begin{eqnarray*}  \nonumber 
 y^0_1=c_2\sqrt{\kappa_2}, \   y^0_\nu=c_2(\sqrt{\kappa_{\nu+1}}-\sqrt{\kappa_\nu})=\frac{c_2({\kappa_{\nu+1}}-{\kappa_\nu})}{2\sqrt{\kappa_\nu^*}}, \\
 \mbox{where} \ \kappa_\nu<\kappa_\nu^*<\kappa_{\nu+1} \  (\nu=2,..., k). \  \mbox{Therefore} \\
 Y(y^0)= \sum\limits_{\nu=1}^{k}y^0_\nu f(\kappa_\nu) =c_2\left[\sqrt{\kappa_2}f(\kappa_1)+ \sum\limits_{\nu=2}^{k} 
 \frac{(\kappa_{\nu+1}-\kappa_\nu)f(\kappa_{\nu})}{2\sqrt{\kappa_\nu^*}}\right]. 
\end{eqnarray*}
Notice that $\kappa_1=a, \ \kappa_{k+1}=1$.
In the limit as  $\theta\to0$  we get
 \begin{eqnarray*}  \nonumber 
\lim_{\theta\to0}Y{(y^0)}=c_2\left[\sqrt{a}f( a)+\frac{1}{2} \int\limits_{a}^{1} 
 \frac{f(\kappa)}{\sqrt\kappa}d\kappa\right]
<c_2\left[\sqrt{a}f( a)+\frac{1}{2} \int\limits_{0}^{1} 
 \frac{f(\kappa)}{\sqrt\kappa}d\kappa\right]<c_3. 
\end{eqnarray*}

\end{document}